\documentclass[12pt]{article}
\usepackage[final]{epsfig}
\usepackage{graphics}
\usepackage{amsmath}
\usepackage{amsfonts}
\usepackage{latexsym}
\usepackage{amssymb}
\usepackage{graphicx}
\usepackage{url}
\usepackage{epstopdf}
\usepackage{xcolor}

\newtheorem{lemma}{Lemma}[section]
\newtheorem{proposition}[lemma]{Proposition}
\newtheorem{remark}[lemma]{Remark}
\newtheorem{example}[lemma]{Example}
\newtheorem{theorem}{Theorem}
\newtheorem{definition}[lemma]{Definition}
\newtheorem{corollary}[lemma]{Corollary}

\begin{document}
\newcommand{\eps}{{\varepsilon}}
\newcommand{\g}{{\gamma}}
\newcommand{\aA}{{\alpha}}
\newcommand{\bB}{{\beta}}
\newcommand{\proofend}{$\Box$\bigskip}
\newcommand{\C}{{\mathbb C}}
\newcommand{\Q}{{\mathbb Q}}
\newcommand{\R}{{\mathbb R}}
\newcommand{\Z}{{\mathbb Z}}
\newcommand{\RP}{{\mathbb {RP}}}
\newcommand{\CP}{{\mathbb {CP}}}
\newcommand{\Tr}{\rm Tr}
\def\proof{\paragraph{Proof.}}

\title{Wire billiards, the first steps}

\author{Misha Bialy\footnote{
School of Mathematical Sciences,
Tel Aviv University, Israel;
bialy@post.tau.ac.il
}
\and
Andrey Mironov\footnote{
Novosibirsk State University and Sobolev Institute of Mathematics,
Novosibirsk, Russia;
mironov@math.nsc.ru
}
\and
Serge Tabachnikov\footnote{
Department of Mathematics,
Penn State University, USA;
tabachni@math.psu.edu}
}

\date{}
\maketitle
\begin{abstract}
	{\it Wire billiard }is defined by a smooth embedded closed curve of non-vanishing curvature $k$ in $\R^n$ (a {wire}). For a class of curves, that we call nice wires, the wire billiard map is area preserving twist map of the cylinder.
In this paper we are investigating whether the basic features of conventional planar billiards extend to this more general situation. 
In particular, we extend Lazutkin's  KAM result, as well as Mather's converse KAM result, to wire billiards. We address the notion of caustics:  for wire billiards, it corresponds to {\it striction curve} of the ruled surface spanned by the chords of the invariant curve.
If the ruled surface is developable this is a genuine caustic.
We found remarkable examples of the wires which are closed orbits of 1-parameter subgroup of $SO(n)$. These wire billiards are totally integrable. Using the theory of interpolating Hamiltonians, we prove that the distribution of impact points of the wire becomes  
uniform with respect to the measure $k^{2/3} dx$ (where $x$ is the arc length parameter), as the length of the chords tends to zero. Applying this result, we prove that the billiard transformation in an ellipsoid commutes with the reparameterized geodesic flow on a confocal ellipsoid: the speed of the foot point of the line tangent to a geodesic equals $k^{-2/3}$, where $k$ is the curvature of the geodesic in the ambient space. We also discuss  perspectives and open problems of this new class of billiards. 
\end{abstract}
\section{Introduction} \label{intro}

In this paper we introduce a billiard-like system that we call {\it wire billiard}.

In the familiar, conventional setting, a billiard table is a domain in Euclidean space $\R^n$ with a smooth boundary $M^{n-1}$. The billiard dynamics describes the motion of a free particle inside the domain with the elastic reflection off the boundary described by the familiar law of geometrical optics: the angle of incidence equals the angle of reflection. In other words, at the moment of collision with $M$, the velocity of the particle is reflected in the tangent hyperplane to $M$ at the collision point.

We propose to replace the hypersurface $M$ by a smooth closed curve $\g$, a {\it wire}. The reflection acts on chords of $\g$: a chord $[\g(x),\g(y)]$ reflects to a chord $[\g(y),\g(z)]$ if they make equal angles with the tangent line to $\g$ at point $y$. More precisely
\begin{equation} \label{eq:law}
\frac{\g(y)-\g(x)}{|\g(y)-\g(x)|} \cdot \dot \g(y) = \frac{\g(z)-\g(y)}{|\g(z)-\g(y)|} \cdot \dot \g(y).
\end{equation}

In general, this reflection law is not deterministic: it defines a correspondence, that is, a multi-valued map. A similar situation was considered in a recent paper \cite{FKM}.

Wire billiards are related to the usual billiards, when the reflection occurs from hypersurfaces.
Consider a smooth curve $\g$ in $\R^n$, and let $M^{n-1}$ be the boundary of  its tubular $\eps$-neighborhood. Consider an instance of the billiard reflection at a point $x \in M$: the incoming ray $\ell$ reflects to the outgoing ray $\ell_1$. Let $v$ and $v_1$ be the unit tangent vectors of the rays and $N$ be the normal to $M$ at point $x$. Then the law of the billiard reflection implies that $v_1 - v$ is proportional to $N$. In the limit $\eps \to \infty$, the point $x$ falls on $\g$, the vector $N$ is orthogonal to $\g$, and the billiard reflection law becomes the law of the wire billiard reflection (\ref{eq:law}).

The literature on mathematical billiards is extensive; see, e.g., the surveys \cite{CM,KT,Tab95,Tab05}. In this paper we make the first steps in the study of wire billiards. In particular, we are investigating whether the basic features of conventional planar billiards extend to this more general situation. 

The contents of the paper are as follows. 

In Section \ref{sec:basic} we examine the generating function and the invariant area form of the wire billiard correspondence. We define a class of curves, called {\it nice}, for which this correspondence is single-valued, and prove, in Theorem \ref{truemap}, that it is an area preserving twist map. Theorem \ref{thm:open} asserts that the niceness condition is open in $C^2$-topology.

In Section \ref{sec:Ma} we prove a version of Mather's theorem: if the curvature of the wire vanishes at some point, then the wire billiard does not have invariant circles (Theorem \ref{sec:Ma}).

Section \ref{sec:Laz} concerns the existence of invariant circles. We show that, similarly to the planar case, if the curvature of the wire $\g$ (not necessarily nice) does not vanish and  $\g$ is smooth enough, then, arbitrarily close to the boundary of the phase cylinder, there exist invariant circles of the wire billiard correspondence. This Theorem \ref{thm: invcircles} is an analog of the celebrated Lazutkin theorem \cite{Laz}.  

Another result of Section \ref{sec:Laz} concerns the limit distribution of the impact points along a wire with non-vanishing curvature $k$ as the length of the chords tends to zero:  Theorem \ref{asympt} asserts that this distribution is uniform with respect to the measure $k^{2/3} dx$ where $x$ is the arc length parameter. This result is  not new, it was obtained in \cite{Gl} by a totally different method; our approach is via the theory of 
interpolating Hamiltonians of Melrose.

Section \ref{sect:caus} concerns analogs of caustics for wire billiards. For usual planar billiards, a caustic is a curve such that if a segment of a billiard trajectory is tangent to it, then the reflected segment is again tangent to this curve. Caustics are the envelopes of the rays that comprise an invariant curve of the billiard transformation of the phase cylinder. 

If a wire billiard has an invariant circle, the respective 1-parameter family of oriented lines forms a ruled surface, but does not necessarily envelope a curve. The striction curve of a ruled surface consists of the points of the rulings closest to the infinitesimally close rulings (if the surface is developable, then this curve is the edge of regression). 
In Theorem \ref{thm:striction} we show that the striction curve of the ruled surface, corresponding to an invariant circle of a wire billiard map, intersects every respective chord of the wire, and not just the lines that they span. This is an analog of the well known fact that the caustics lie inside the billiard table. 

 Section \ref{sec:tot} is devoted to constructions of non-planar totally integrable wire billiards whose phase cylinder is foliated by invariant circles; they are analogs of circles in the plane. Theorem \ref{thm:niceint} provides such an example, a  nice curve in $\C^2 = \R^4$ given by the formula
 $$
 \g(t) = (e^{it}, \eps e^{imt}),\ m \ge 2
 $$
 with sufficiently small $\eps$. This curve is an orbit of a 1-parameter subgroup of $SO(4)$, and the invariant circles consist of the chords $[\g(t) \g(t+d)]$ for a fixed $d$. 
 
 Theorem \ref{thm:caustic} asserts that to such an invariant circle there corresponds a genuine caustic (that is, the ruled surface is developable) if and only if the relation
 $$
 \tan\left(m\frac{d}{2}\right)=m\tan \left(\frac{d}{2}\right)
 $$
 holds. This ubiquitous  equation appeared in the study of billiards and other topics before, see Remark \ref{rm:Gut}. 
 
In Section \ref{sect:ell} we apply Theorem \ref{asympt} to billiards inside ellipsoids in $\R^n$. These billiards are completely integrable: the $(2n-2)$-dimensional symplectic space of oriented lines is foliated by invariant Lagrangian submanifolds, consisting of the lines tangent to fixed $n-1$ confocal quadrics. As a result, the billiard transformations on the rays, that intersect two two confocal ellipsoids, commute; see, e.g. \cite{DR,Tab95} for this material.
 
 The geodesic flow on an ellipsoid can be considered as a flow on the lines that are tangent to this ellipsoid. This flow is also completely integrable, with the same integrals as the billiard transformation. Theorem  \ref{thm:repar} asserts that the billiard transformation in an ellipsoid commutes with the reparameterized geodesic flow on a confocal ellipsoid: the speed of the foot point of the line, tangent to a geodesic, equals $k^{-2/3}$, where $k$ is the curvature of this geodesic.

The final Section \ref{sec:problems} contains a number of open problems and conjectures.

\bigskip
{\bf Acknowledgements}. MB was supported by ISF grant 162/15, AM was supported by RSF grant 19-11-0004 , ST was supported by NSF grant DMS-1510055.


\section{Wire billiards: generating function, twist condition, invariant area form} \label{sec:basic}

Let $\g$ be a curve in $\R^n$. Chose the arc-length parametrization of $\gamma$. The {\it wire billiard} on $\g$ is defined by its generating function
$L(x,y)=|\g (y)-\g (x)|$, by the usual formula $T: (\g (x),\g (y)) \mapsto (\g (y),\g (z))$ if

\begin{equation}\label{genfunct}
L_2(x,y) + L_1(y,z) = 0.
\end{equation}
Here and below we shall write subindices 1 and 2 for partial derivatives of the function $L$ with respect to first and the second argument respectively.

Defined in this way, $T$  is not necessarily a map, but rather a relation (a multivalued map). One may consider $T$ as a non-deterministic billiard-like system, cf. \cite{FKM}.
Below we specify conditions on $\g$ that guarantee that $T$ is a map.

Let $\aA$ and $\bB$ be the angles made by a chord $[\g(x), \g(y)]$ with $\g$ at $\g(x)$ and $\g(y)$, respectively. Let $\omega$ be the 2-form $\sin \alpha\ d\alpha \wedge dx$.

\begin{lemma}  \label{symp1}
The relation $T$ is symplectic: $T^*(\omega)=\omega$.
\end{lemma}

\proof
It is well known that the 2-form
$$
L_{12}(x,y) dx \wedge dy
$$
is invariant under the discrete Lagrangian dynamics, see, e.g., \cite{Ve}. For completeness, here is the proof.

Take the exterior derivative of (\ref{genfunct})
$$
L_{12}(x,y)dx + L_{22}(x,y) dy +
L_{11}(y,z) dy + L_{12}(y,z) dz = 0
$$
and wedge multiply by $dy$ to obtain
$$
L_{12}(x,y) dx \wedge dy = L_{12}(y,z) dy \wedge dz,
$$
as needed.

Since,
$$
L^2=(\gamma(y)-\gamma(x))\cdot (\gamma(y)-\gamma(x)),
$$
we have the first derivatives of $L$:
$$
	L_1(x,y)=-\frac{\dot\gamma(x)\cdot(\gamma(y)-\gamma(x))}{L} = -\cos\aA,\
$$
\begin{equation}\label{first}
	L_2(x,y)=\frac{\dot\gamma(y)\cdot(\gamma(y)-\gamma(x))}{L}=\cos\bB.
	\end{equation}
Hence
$$
L_{12}(x,y)= \sin \aA \frac{\partial \aA}{\partial y},
$$
and therefore
$$
L_{12} dx \wedge dy = \sin \aA \frac{\partial \aA}{\partial y} dx \wedge dy = - \omega,
$$
as claimed.
\proofend

The space of oriented lines in $\R^n$ has a symplectic structure $\Omega = dp \wedge dq$, where $q$ is the unit vector along the line and $p$ is its closest point to the origin. The notation $dp \wedge dq$ means $\sum dp_i \wedge dq_i$.

In fact, the form $\Omega$ does not change if, instead $p$ one takes another point on the line. Indeed, if $p_1=p+f(p,q) q$ then
$$
dp_1 \wedge dq = dp \wedge dq + df \wedge q dq = dp \wedge dq,
$$
the latter equality due to the fact that $q dq=0$ since $q\cdot q =1$.

Let $\mathcal{C}$ be the phase space of the wire billiard, the space of chords of the curve $\g$, that is, of oriented lines that intersect $\g$ twice.

\begin{lemma} \label{symp2}
The restriction of $\Omega$ to $\mathcal{C}$ is $\omega$.
\end{lemma}

\proof Let $q$ be the unit vector along a chord and let $p= \g(x)$. Then $\cos \aA = \dot\g(x) \cdot q$. One has
$$
\omega = dx \wedge d (\cos \aA) = dx \wedge (\dot\g(x) dq) = dp \wedge dq,
$$
as claimed.
\proofend

Given $x$ and $y$, consider the 2-planes $\pi_{xy}=Span\{\dot\gamma(x),\gamma(y)-\gamma(x)\}$ and $\pi_{yx}=Span\{\dot\gamma(y),\gamma(x)-\gamma(y)\}$. These planes share a line and span 3-dimensional space $V$. Let $\varphi (x,y)$ be the angle between these two planes.

\begin{lemma} \label{angle}
One has
$$
L_{12}(x,y) = \frac{\cos\varphi \sin \aA \sin \bB}{L}.
$$
\end{lemma}

\proof
Using formulas (\ref{first}), let us compute the mixed  second partial derivative:
$$
\frac{\partial}{\partial y}\left(L_1 \right)=-\frac{\dot\gamma(x)\cdot\dot\gamma(y)}{L}+\frac{\dot\gamma(x)\cdot(\gamma(y)-\gamma(x))}{L^2}L_2=
$$
$$
=-\frac{\dot\gamma(x)\cdot\dot\gamma(y)}{L}+\frac{\dot\gamma(x)\cdot(\gamma(y)-\gamma(x))}{L^2}\cdot
\frac{\dot\gamma(y)\cdot(\gamma(y)-\gamma(x))}{L}=
$$
$$
=\frac{-\left(\dot\gamma(x)\cdot\dot\gamma(y)\right)(\gamma(y)-\gamma(x))^2+(\dot\gamma(x)\cdot(\gamma(y)-\gamma(x)))(\dot\gamma(y)\cdot(\gamma(y)-\gamma(x)))}{L^3}.
$$
Note that the numerator of the last expression can be computed in the 3-space $V$:
$$
[\dot\gamma(x) \times (\gamma(y)-\gamma(x))]\cdot [ \dot\gamma(y)\times (\gamma(x)-\gamma(y))],
$$
due to the identity that holds for any four vectors in $\R^3$:
$$
(a\times b) \cdot (c\times d) = (a\cdot c) (b\cdot d) - (a\cdot d) (b\cdot c).
$$
Note also that
$$
[\dot\gamma(x) \times (\gamma(y)-\gamma(x))]\cdot [ \dot\gamma(y)\times (\gamma(x)-\gamma(y))] = L^2 \cos\varphi \sin \aA \sin \bB.
$$
This implies the statement of the lemma.
\proofend

\begin{remark}
{\rm If the curve $\g$ is planar,  the formula of Lemma \ref{angle}, with $\varphi=0$, is well known in the literature on billiards.
For space curves, this formula, in its relation with the symplectic form on the space of oriented lines, appeared in \cite{Pohl}.
}
\end{remark}

We can now formulate sufficient condition for the wire billiard relation to be a mapping.

\begin{definition} \label{nice}
{\rm
A a smooth closed curve $\g$ is called {\it nice} if it satisfies the following conditions:
\begin{enumerate}
\item Any line intersects $\g$ in at most two points, and if so, the line makes non-zero angles with the curve at both intersection points;
\item The curvature $k$ of $\g$ is strictly positive;
\item For any $x,y$, the 2-planes $\pi_{xy}$ and $\pi_{yx}$ are not orthogonal.
\end{enumerate}
}
\end{definition}

For example, a planar strictly convex (i.e., with positive curvature) curve is nice.
Note that a nice curve in $\R^3$ is unknotted: the chords connecting a fixed point with other points constitute a spanning disc of the curve.

\begin{theorem} \label{truemap}
For a nice curve,  $T$ is an area preserving twist map of the phase cylinder $\mathcal{C}$.
\end{theorem}

\proof
Lemmas \ref{symp1} and \ref{angle} imply that the phase area form $\omega$ is non-vanishing.

Lemma \ref{angle} implies that the twist condition
$$
L_{12}(x,y) >0
$$
holds. Indeed,
$$
L_{12}(x,y)   \neq 0,
$$
and the sign must be positive because the angle $\varphi (x,y)$ is small when $x$ and $y$ are close to each other.

We have
$$
0 <L_{12}(x,y)  = \frac{-\partial (\cos \aA)}{\partial y}= \sin \aA \frac{\partial \aA}{\partial y}.
$$
Therefore, for every fixed $x$, the angle $\aA$ increases monotonically as $y$ increases. This implies that, given a chord $[\g(x) , \g(y)]$, there is a unique $z$ such that the chord $[\g(y), \g(z)]$ makes the same angle with the curve at point $\g(y)$ as the chord $[\g(x), \g(y)]$. That is, $T$ is a (single-valued) map.
\proofend

\begin{remark}
{\rm Suppose that  $\g \subset \R^n$ is a nice curve.
The total phase area of $\mathcal{C}$ equals $2 |\g|$, where $|\g|$ is the length of the curve $\g$. Just like in the plane, this follows from integrating the form $\omega = \sin \alpha\ d\alpha \wedge dx$.

If $\g$ is a convex planar curve, the integral of the generating function $L$ over the phase space with respect to the area form $\omega$ equals $2\pi A$, where $A$ is the area bounded by $\g$, see \cite{Tab05}.

In the 3-dimensional case, an analog of this result was discovered in \cite{Pohl}. For a non-oriented line $\ell$ in $\R^3$, let $\lambda(\ell)$ be the absolute value of the linking number of $\ell$ with $\g$. Then
$$
\int_{\g \times \g} L(x,y)\ \omega = 2 \int \lambda^2(\ell)\ \Omega\wedge \Omega,
$$
where $\Omega\wedge \Omega$ is the symplectic volume form on the space of non-oriented lines in $\R^3$. For simple curves of a fixed length, this quantity is maximal for planar circles, see \cite{BP}.
}
\end{remark}

The next result shows that niceness is an open property.

\begin{theorem} \label{thm:open}
A sufficiently small perturbation in the $C^2$-topology of a smooth nice curve is a nice curve.
\end{theorem}

\proof
Let $\g$ be a nice curve.
The curvature of $\g$ is a continuous function on the curve, and a small perturbation keeps the curvature positive.

Compactify the phase cylinder ${\mathcal C}$, that is, the space of chords $\g \times \g - {\rm Diag}$, by replacing the diagonal with the unit tangent bundle to $\g$, that is, two circles consisting of unit tangent vectors to $\g$, in both directions.

The angle $\varphi$ between the planes $\pi_{xy}$ and $\pi_{yx}$ extends continuously to this compactification: it takes zero values on the added circles. Indeed, as $y\to x$, the limiting position of the planes $\pi_{xy}$ and $\pi_{yx}$ is the osculating plane of the curve $\g$ at point $\g(x)$.
Since $\varphi < \pi/2$, a small perturbation of $\g$ preserves this inequality.

It remains to show that the twist condition
$$
L_{12}(x,y)  = \frac{\cos\varphi \sin \aA \sin \bB}{L} >0
$$
also survives under a small perturbation. It is suffices to show that, for a sufficiently small perturbation of $\g$,
\begin{equation} \label{strongineq}
\frac{\sin \aA \sin \bB}{L^2} >0.
\end{equation}

We argue that the functions
$$
\frac{\sin \aA }{L} \ \ {\rm and}\ \ \frac{ \sin \bB}{L}
$$
extend continuously to the compactification of ${\mathcal C}$. Consider the first function; the second one is dealt with similarly.

One has
$$
\dot \g \cdot \dot \g =1,\ \ddot \g \cdot \ddot \g = k^2.
$$
Differentiating, one gets
\begin{equation} \label{manydots}
\dot \g \cdot \ddot \g = 0,\ \dot \g \cdot \dddot \g = -k^2,
\ \ddot \g \cdot \dddot \g = \dot k k,\ \dot \g \cdot \ddddot \g = -3\dot k k.
\end{equation}

One  has
$$
\cos \aA = \frac{\dot \g(x) \cdot (\g(x+\eps)-\g(x))}{|\g(x+\eps)-\g(x)|}.
$$
Considering Taylor expansion modulo $\eps^4$ and using  formulas (\ref{manydots}), we obtain
$$
\cos \aA = 1 - \frac{\eps^2}{8} k^2 - \frac{\eps^3}{12} \dot k k + (\eps^4),
$$
and hence
\begin{equation} \label{eq:angle}
\aA = \frac{\eps}{2} k + \frac{\eps^2}{6}\dot k + (\eps^3).
\end{equation}
This implies
$$
\lim_{\eps \to 0} \frac{\sin \aA}{L} = \frac{k}{2}.
$$

Therefore $\sin \aA /{L}$ is a positive function on a compact space, and it remains positive under a sufficiently small perturbation.
\proofend

\begin{corollary} \label{cor:open}
A sufficiently small perturbation in the $C^2$-topology
of any planar smooth strictly convex curve is nice.
\end{corollary}

\begin{example} \label{ex:perp}
{\rm Let $\g_1 = (\cos t, \sin t, 0, 0)$ and $\g_2 = (0, 0, \cos s, \sin s)$ be two curves in $\R^4$. Then every chord $[\g_1(t), \g_2(s)]$ is orthogonal to the curves $\g_1$ and $\g_2$ at its endpoints. If $\g$ is a curve in $\R^4$ that contains pieces of $\g_1$ and $\g_2$, then $\g$ has an open 2-dimensional set of 2-periodic wire billiard trajectories. Of course, in this case the wire billiard correspondence is not a single-valued map; it is rather infinitely-valued one.
}
\end{example}

\section{Glancing wire billiards} \label{sec:Ma}

The extension of J. Mather theorem \cite{Ma}  for wire billiards is as follows. As usual, this theorem implies the existence of glancing orbits on the phase cylinder $\mathcal{C}$, that is, the orbits that make arbitrarily small angles with the tangent vectors to $\g$, both positive and negative.

  \begin{theorem} \label{thm:Ma}
  	Let $k$ be the curvature of the curve $\gamma$ in $\R^n$. If $k$ vanishes at one point then the billiard ball map of the phase cylinder $\mathcal{C}$ has no non-contractible invariant curves.
  \end{theorem}

\proof
Differentiating (\ref{first}), we compute the other second derivatives of $L$:
\begin{equation*}
\begin{split}
L_{11}(x,y)&=-\frac{\ddot\gamma(x)\cdot(\gamma(y)-\gamma(x))-1}{L}+
\frac{\dot\gamma(x)\cdot(\gamma(y)-\gamma(x))}{L^2}L_1\\
&=\frac{1-\ddot\gamma(x)\cdot(\gamma(y)-\gamma(x))}{L}-\frac{\cos^2\alpha}{L}.
\end{split}
\end{equation*}
Hence
\begin{equation}\label{secondx}
L_{11}(x,y)=\frac{\sin^2\alpha}{L}-k(x)\left( n(x)\cdot\frac{\gamma(y)-\gamma(x)}{L}\right),
\end{equation}
where $n(x)$ is the unit normal to $\g$ at $\g(x)$.
Similarly,
\begin{equation}\label{secondy}
L_{22}(x,y)
=\frac{\sin^2\beta}{L}-k(y)\left( n(y)\cdot\frac{\gamma(x)-\gamma(y)}{L}\right).
\end{equation}

Using these formulas, the proof proceeds analogously to the planar case.
If there is a non-contractible invariant curve on $\mathcal{C}$ then, by Birkhoff's theorem for twist maps, it must be a graph of a Lipshitz function. Denote by $f:\gamma\rightarrow\gamma$ the induced homeomorphism
on the invariant curve. Next we write equation (\ref{genfunct}) for the piece of the orbit of point $q\in\gamma$, namely, $(f^{-1}(q),q,f(q))$:
$$
L_2(f^{-1}(q),q) + L_1(q,f(q)) = 0.
$$
Differentiating with respect to $q$, we have:

$$
L_{11}(q,f(q))+L_{22}(f^{-1}(q),q)=
-\left[L_{12} (q,f(q))f^{'}(q)+L_{12} (f^{-1}(q),q)(f^{-1})^{'}(q)\right].
$$
Note that, by the twist condition, the right hand side is negative. On the other hand,
choose $q\in\gamma$ to be the point of zero curvature. Then it follows from (\ref{secondx}) and (\ref{secondy}) that the left hand side is positive. This contradiction proves the theorem.
\proofend

\section{Lazutkin parameter} \label{sec:Laz}

Let $\g(x)$ be arc length parameterized curve in $\R^n$, $k(x)>0$ be its curvature, $x_{\pm} = x\pm\eps$.
Let $\g_{\pm}$ and $k_{\pm}$ mean $\g(x_{\pm})$ and $k(x_{\pm})$ respectively. Let $\aA_{\pm}$ be the angle made by the curve at points $\g_{\pm}$ with the chord $[\g({x_-}),\g({x_+})]$.

The Lazutkin coordinates $(u,v)$ in the phase space of the wire billiard are given by the equations
 $$
 du=k^{2/3} dx,\ v= k^{-1/3} \sin\left(\frac{\aA}{2}\right).
 $$
 Clearly, up to a factor, one has $\omega=v dv \wedge du$.

 Let $T(u,v)=(u_1,v_1)$.


 \begin{proposition} \label{approx}
 One has:
 $$
 u_1 = u + 4 v\ {\rm mod}\ v^3,\  v_1 = v\ {\rm mod}\ v^4.
 $$
 \end{proposition}

 \proof
 Let us prove the first equality. One has
 $$
 u_1 - u = \int_0^\eps k(x+\tau)^{2/3} d \tau.
 $$
Since
$$
k(x+\tau)^2 = \ddot\g(x+\tau) \cdot \ddot\g(x+\tau) = k^2 + 2\tau \dot k k + (\tau^2),
$$
where $k$ is shorthand for $k(x)$, we have
$$
k(x+\tau)^{2/3} = k^{2/3} \left( 1 + \frac{2}{3} \frac{\dot k}{k} \tau \right) + (\tau^2).
$$
Integrating with respect to $\tau$  from 0 to $\eps$, we obtain
$$
u_1 - u = \eps k^{2/3} + \frac{\eps^2}{3} \dot k k^{-1/3} + (\eps^3).
$$

 On the other hand, we know from (\ref{eq:angle}) that
 $$
\sin\left(\frac{\aA}{2}\right) = \frac{\eps}{4} k + \frac{\eps^2}{12} \dot k + (\eps^3),
 $$
 hence
 $$
 4v = 4k^{-1/3} \sin\left(\frac{\aA}{2}\right) = u_1 - u\ {\rm mod}\ (\eps^3),
 $$
 as needed.

To prove the second equality, we show that $v_+ = v_-\ {\rm mod}\ \eps^4$.

Using formulas (\ref{manydots}), we calculate:
 $$
 \g_+-\g_- = 2\eps \left(\dot\g +\frac{\eps^2}{6}\dddot\g +(\eps^4)\right),
 $$
 and
 $$
 |\g_+-\g_-| = 2\eps \left(1-\frac{\eps^2}{6}k^2+(\eps^4)\right),
 $$
hence
 $$
 \frac{\g_+-\g_-}{|\g_+-\g_-|} = \dot\g + \frac{\eps^2}{6} (k^2\dot\g + \dddot\g) + (\eps^4).
 $$

 Next,
 $$
 \dot\g_\pm = \dot\g\pm \eps \ddot\g+\frac{\eps^2}{2}\dddot\g\pm \frac{\eps^3}{6}\ddddot\g+(\eps^4),
$$
hence
$$
\cos\aA_\pm = \dot\g_\pm\cdot \frac{\g_+-\g_-}{|\g_+-\g_-|} = 1 - \frac{\eps^2}{2}k^2 \mp \frac{\eps^3}{3}\dot k k +(\eps^4),
$$
and thererefore
$$
\aA_\pm = |\eps| \left(k \pm \frac{\eps}{3} {\dot k}\right) + (\eps^3).
$$
Thus
$$
\sin \left(\frac{\aA_\pm}{2}\right) = \frac{|\eps|}{2} \left(k\pm\frac{\eps}{3}\dot k +\eps^2 f(x) \right) +(\eps^4),
$$
where $f(x)$ is some function.

Similarly,
$$
k_\pm^{-1/3} = k^{-1/3} \left(1 \mp \frac{\eps}{3} \frac{\dot k}{k} + \eps^2 g(x) \right) + (\eps^3).
$$
It follows that
$$
k_\pm^{-1/3} \sin \left(\frac{\aA_\pm}{2}\right) = \frac{|\eps|}{2} k^{2/3}
\left[ 1 + \eps^2 \left(\frac{ f}{k} + g - \frac{1}{9} \frac{\dot k^2}{k^2}\right) \right] +
(\eps^4),
$$
hence
$$
k_+^{-1/3} \sin \left(\frac{\aA_+}{2}\right) = k_-^{-1/3} \sin \left(\frac{\aA_-}{2}\right)\ {\rm mod}\ \eps^4,
$$
as needed.
 \proofend

Proposition \ref{approx} says that, near the boundary of the phase cylinder, the map is close to an integrable one. It
has important consequences.

The first one is the existence of invariant circles of the wire billiard correspondence. Assume that $\g$ is an infinitely  smooth closed curve with non-vanishing curvature, but not necessarily nice.

\begin{theorem} \label{thm: invcircles}
The wire billiard correspondence has invariant circles arbitrarily close to the boundary of the phase space ${\mathcal C}$. The union of these invariant circles has a positive measure.
\end{theorem}

\proof
The arguments of Theorem \ref{thm:open} and Proposition \ref{approx} apply, showing that, for sufficiently small angles, the wire billiard map is still a well defined area preserving twist map, close to an integrable one.

Now one applies a KAM theory result, the Lazutkin theorem (Theorem 2 of \cite{Laz}), that states that an area preserving twist map
$$
u_1 = u + v + f(u,v)  v^m,\  v_1 = v + g(u,v) v^{m+1}
$$
with $m \ge 1$ and sufficiently smooth functions $f,g$ has invariant curves arbitrarily close to the boundary $v=0$.
By Proposition \ref{approx}, $m=3$, as in the familiar situation of planar billiards.
\proofend


 Secondly, as in \cite{MM}, Proposition \ref{approx} implies the limit distribution of the impact points of gliding billiard trajectories.

\begin{theorem} \label{asympt}
Let $\g$ be an arc of a curve. Consider a wire billiard trajectory on $\g$ with $n$ impact points whose endpoints are the endpoints of the arc $\g$ (that is, the longest polygon, properly inscribed in $\g$). Let $\ell$ be the length of the arc $\g$ and $\ell_n$ be the length of the billiard trajectory. Then\\
1) In the limit $n\to\infty$, the distribution of the impact points is uniform with respect to the measure $k^{2/3} dx$. \\
2) One has
$$
\lim_{n\to\infty} n^2 (\ell - \ell_n) = \frac{1}{24} \int_{\g} k^{2/3} dx.
$$
\end{theorem}

The proof follows from the fact that $v^2/2$ is an interpolating Hamiltonian, and the flow of its Hamiltonian vector field $\partial/\partial u$ approximates the billiard map as $v \to 0$, see \cite{MM,Me,PS}.

The asymptotic formula of Theorem \ref{asympt} is known:
it was stated by Fejes T\'oth \cite{FT} and proved by Gleason \cite{Gl}; see also  \cite{Gr,McV}. These proofs do not use symplectic geometry and the theory of interpolating Hamiltonians.

\begin{remark}
{\rm A more invariant description of the Lazutkin parameter is as follows. Consider a curve $\g$ with non-vanishing curvature, and introduce a cubic form on $\g$ as follows.

Let $v$ be a tangent vector to $\g$ at point $x \in \g$. Include $v$ into a tangent vector field on $\g$ and flow along this field for time $\eps$ to point $y$. Consider the difference of the length of the arc $xy$ and the chord $[x,y]$. It depends cubically on  $\eps$. Divide by $\eps^3$ and take limit $\eps \to 0$ to obtain a function of $x$ and $v$ which is cubic in $v$, that is, a cubic form on $\g$. A calculation yields that this form is $k^2 dx^3$. The cube root is Lazutkin's parameter $du = k^{2/3} dx$.
}
\end{remark}

\begin{remark}
{\rm The theory of interpolating Hamiltonians provides an asymptotic expansion
$$
\ell_n \sim l + \sum_{j=1}^{\infty} \frac{c_j}{n^{2j}}.
$$
Theorem \ref{asympt} provides the value of the coefficient $c_1$. What is the value of $c_2$? In the planar case, it is given in  \cite{MM}.
}
\end{remark}

Now consider the billiard inside a closed smooth strictly convex hypersurface $M \subset \R^n$. The billiard trajectories making small angle with $M$ tend to the geodesics on $M$, and the distribution of the impact points tend to the uniform with respect to the measure $k^{2/3} dx$ on the respective geodesic. See \cite{GP,KP,Po} for billiard flow near the smooth boundary of a strictly convex billiard table.

\section{Invariant curves and caustics} \label{sect:caus}

An invariant circle of the wire billiard map is a one-parameter family of chords of the curve $\g$ that is invariant under the map $T$. In the conventional planar situation, this one-parameter family of lines has an envelope, a billiard caustic. In the multi-dimensional case, a general one-parameter family of lines generates a ruled surface. If the lines have an envelope, that is, are tangent to a spacial curve, this surface is developable. See, e.g., \cite{Lan} for the geometry of ruled surfaces.

Let $\ell(t)$ be a one-parameter family of lines in $\R^n$. Consider two lines, $\ell(t)$ and $\ell(t+\eps)$, and assume that they are not parallel. Then the lines have a unique common perpendicular, say $[x_0 y_0]$, where $x_0 (t,\eps) \in \ell(t), y_0(t,\eps) \in \ell(t+\eps)$. Let $x(t) \in \ell(t)$ be the limiting position of $x_0 (t,\eps)$ as $\eps \to 0$. The curve $x(t)$  is called {\it the striction curve} of the ruled surface. If the family has an envelope, then its striction curve is this envelope.

If the one-parameter family of lines corresponds to an invariant curve of wire billiard, it is natural to consider the striction curve of respective ruled surface to be an analog of billiard caustic.
If the ruled surface is developable, we call the envelope of the lines
a wire billiard caustic.

\begin{figure}[hbtp]
\centering
\includegraphics[height=1.5in]{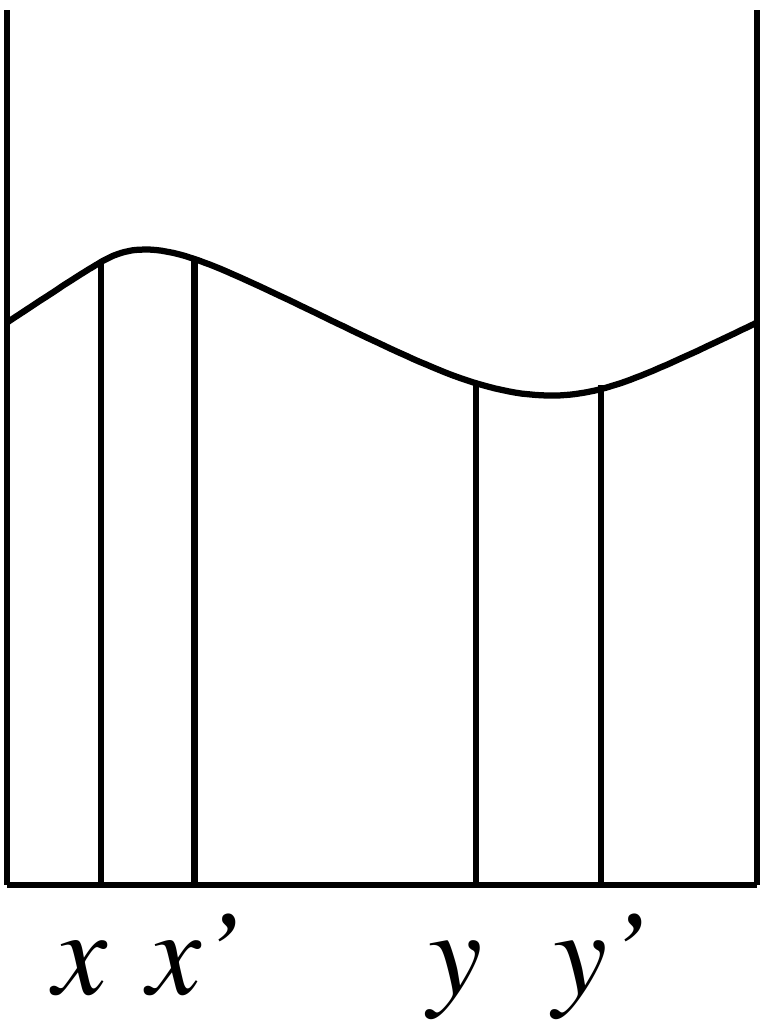} \qquad\qquad
\includegraphics[height=1.5in]{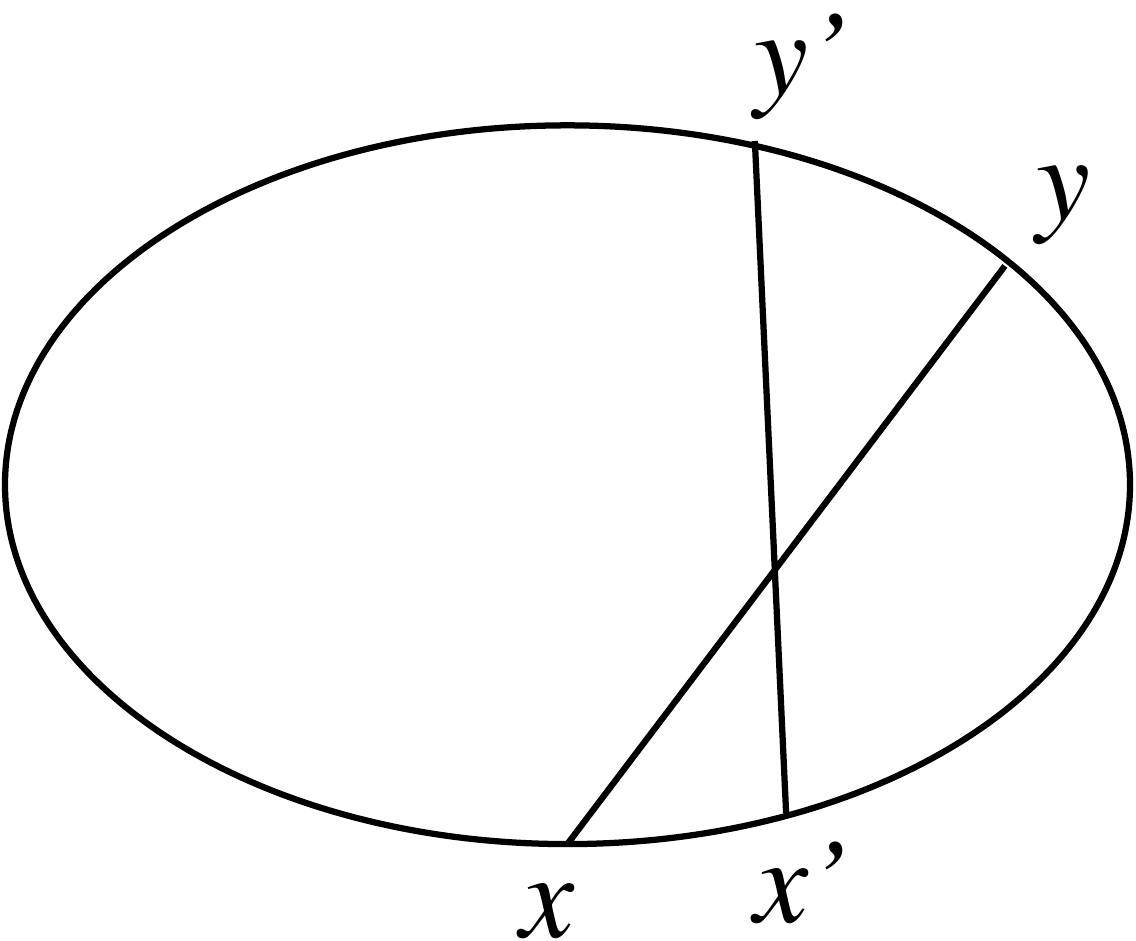}
\caption{Left: an invariant circle of the billiard map with $x=x(t), x'=x(t+\eps), y=y(t), y'=y(t+\eps)$. Right: the respective chords intersect.}
\label{intersect}
\end{figure}

It is well known that the caustics of convex planar billiards lie inside the billiard table. This follows from the above-mentioned theorem of Birkhoff: an invariant circle of the billiard map is a graph. Indeed, if $[x(t) y(t)]$ is a one-parameter family of chords corresponding to an invariant circle then, for small $\eps$, the points $x(t), x(t+\eps), y(t), y(t+\eps)$ lie on the billiard curve in this cyclic order. Therefore the lines   $(x(t) y(t))$ and $(x(t+\eps) y(t+\eps))$ intersect inside the billiard table, see Fig. \ref{intersect}.

A similar property holds for wire billiards.

Let $\g \subset \R^n$ be a nice curve and let the wire billiard map have an invariant curve. Let $[\g(t), \g(f(t))]$ be the respective family of chords of $\g$, where $f:\g\rightarrow\g$ is the corresponding diffeomorphism of $\g$, and let $\delta(t)$ be the striction curve. Write $R(t)= \g(f(t))-\g(t)$.
Then the ruled surface $S$ corresponding to the invariant curve is parameterized as $S(t,s)=\g(t)+sR(t),\ s\in\R$.

\begin{theorem} \label{thm:striction}
1) The ruled surface $S$ corresponding to the invariant curve is not cylindrical:
	$$
	R\wedge \dot R \neq 0;
	$$
	2)
The striction curve $\delta$ intersects all the chords $[\g(t), \g(f(t))]$ (and not merely the lines that they generate).
\end{theorem}

\proof
Since
$$
R(t)= \g(f(t))-\g(t),
$$
we have
$$
\dot {R}(t)=\dot{\g}(f(t))\dot {f} (t)-\dot{\g }(t).
$$

Consider the 3-dimensional space $V$ spanned by the vectors  $R, \dot{\g }(t)$, and $\dot{\g}(f(t))$. It suffices to show that the vectors $R$ and $\dot R$ are linearly independent in $V$.

We have:
\begin{equation}\label{RdotR}
R\times \dot R =[\dot\g(t)\times R]+\dot f(t)[R \times \dot\g(f(t))].
\end{equation}
The two vector products in the last expression are the oriented normals to the planes $\pi_{xy}$ and $\pi_{yx}$, where $x=\g(t),y=\g(f(t))$. By definition of $f$, the derivative $\dot f$ is positive, and we conclude that the angle between the summands is acute. Thus $R\times \dot R\neq 0$, and the two vectors are linearly independent.

In order to prove the second claim, we again work in the space $V$.
Recall the formula for the striction point:
$$
\delta(t)=\g(t)+s^*R(t), \quad s^*=\frac{	[R\times \dot R]\cdot [\dot\g\times R]}{|[R\times \dot R]|^2}.
$$
We want to show that $0<s^*<1$.

Using (\ref{RdotR}), we obtain
$$
[R\times \dot R]\cdot [\dot\g(t)\times R]= |\dot\g(t)\times R|^2 + \dot f(t) [R \times \dot\g(f(t))]\cdot [\dot\g(t)\times R] >0
$$
since the normals to the planes $\pi_{xy}$ and $\pi_{yx}$ make an acute angle. Thus $s^* > 0$.

Finally, we need to show that $s^*<1$, that is, that $|\dot\g\times R|<|R\times \dot R|.$ This again
follows  from (\ref{RdotR}) since the angle between the summands in this formula is acute.
\proofend

Let $\delta_1(t)$ be a (germ of a) curve and $v_1(t)$ be a unit vector field along it, and let $\delta_2(t)$ and $v_2(t)$ be another such pair. Then we have two one-parameter families of oriented lines $\delta_i(t)+sv_i(t),\ i=1,2$. Assume that the first family reflects to the second one in a curve $\g$, and the reflection occurs at point $\g(t)$. Let $a_1(t)$ and $a_2(t)$ be the distances from $\delta_1(t)$ and $\delta_2(t)$ to $\g(t)$.

\begin{proposition} \label{prop:string}
One has the identity:
\begin{equation} \label{eq:stringtion}
(a_1 + a_2)' = \delta_2' \cdot v_2 - \delta_1' \cdot v_1,
\end{equation}
where prime is the derivative with respect to $t$.
\end{proposition}

\proof
One has
\begin{equation} \label{eq1}
\g(t) = \delta_1(t) + a_1(t) v_1(t) = \delta_2(t) - a_2(t) v_2(t),
\end{equation}
and (the reflection law)
\begin{equation} \label{eq2}
\g'(t) \cdot (v_2(t) - v_1(t)) =0.
\end{equation}
Equation (\ref{eq1}) implies
$$
\g' = \delta_1' + a_1' v_1 + a_1 v_1'\ \ {\rm and} \ \ \g' \cdot v_1 = \delta_1' \cdot v_1 + a_1',
$$
where we used that $v_1' \cdot v_1 =0$ since $v_1$ is unit. Likewise,
$$
\g' \cdot v_2= \delta_2' \cdot v_2 - a_2'.
$$
Then equation (\ref{eq2}) implies that
$$
(a_1 + a_2)' = \delta_2' \cdot v_2 - \delta_1' \cdot v_1,
$$
as claimed.
\proofend

In particular, one may take as $\delta$ the striction curve of the respective ruled surface.

\begin{figure}[hbtp]
\centering
\includegraphics[height=1.5in]{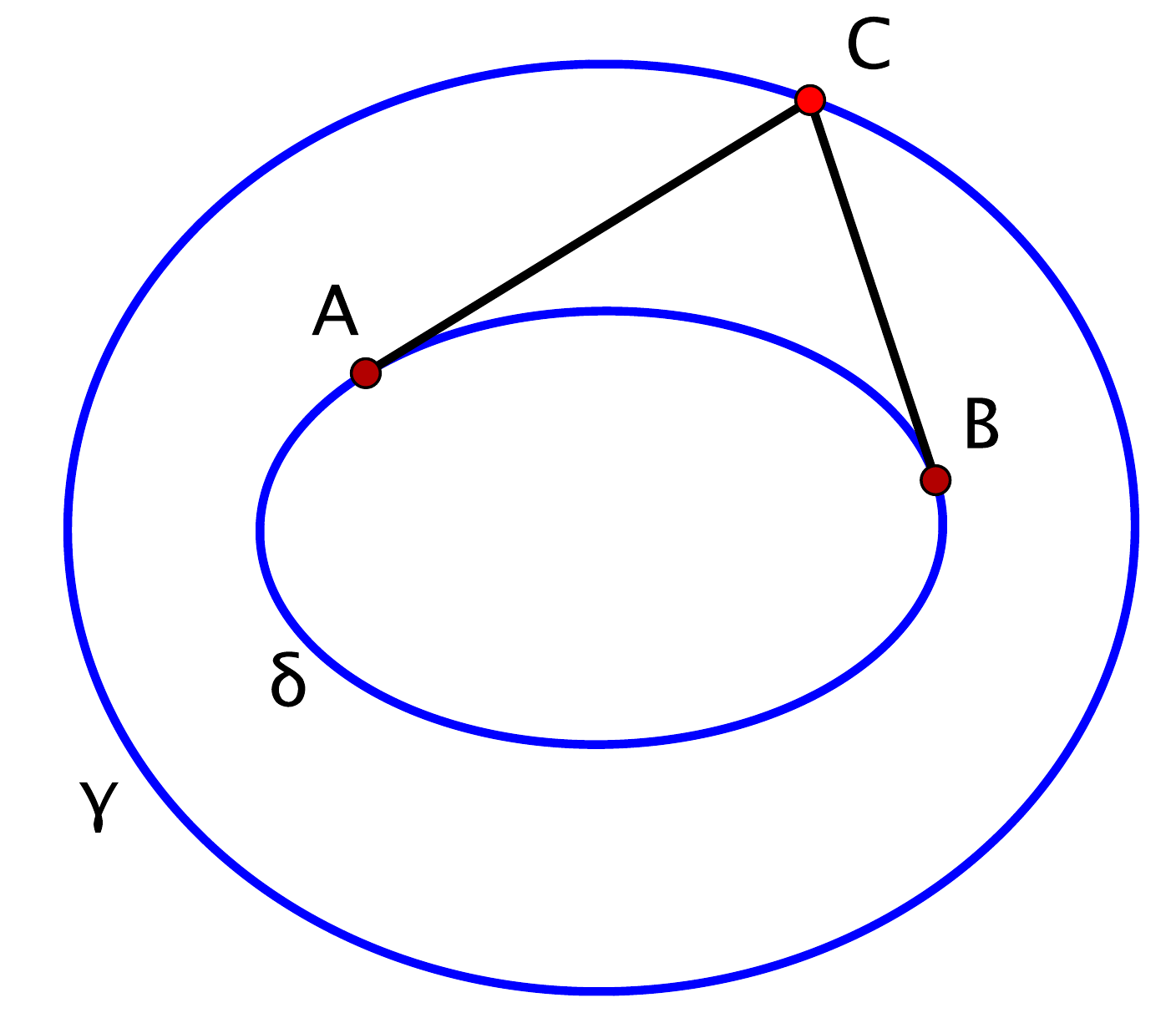}
\caption{The quantity $|AC| + |BC| - |\stackrel{\smile}{AB}|$ is constant along the caustic $\delta$: indeed,
$\delta_2' \cdot v_2 - \delta_1' \cdot v_1$ is the rate of change of the  length of the arc.}
\label{string}
\end{figure}

If the vector field $v_i(t)$ is tangent to the curve $\delta_i$, $i=1,2$, and the curves $\delta_1$ and $\delta_2$ are parts of the same caustic,
then equation (\ref{eq:stringtion}) gives the quantity, conserved along the invariant circle and familiar from the string construction in the planar case, see Fig. \ref{string}.

On the other hand, if $v_i$ is orthogonal to $\delta_i$ for $i=1,2$, then the sum of lengths $a_1+a_2$ is constant.

\section {Totally integrable wire billiard} \label{sec:tot}

{In this Section we work out an important example of totally integrable wire billiard. Following A.Knauf \cite{knauf}, by total integrability we mean the case when the phase space cylinder is foliated by invariant curves which are graphs.}

Consider a closed curve $\gamma\subset\R^n$ which is an orbit of one-parametric subgroup of the group $SO(n)$:
$$
\g(t)=e^{At}\g_0, \ {\rm where} \  A\in so(n),\ 
$$
is a skew-symmetric matrix.
In this case
$$
\dot\g(t)=Ae^{At}\g_0, \ddot\g(t)=A^2e^{At}\g_0,
$$
hence, for almost all $\g_0$, the curve $\gamma$ is a smoothly embedded with non-vanishing curvature in $\R^n$.
\begin{example}Let
$\gamma$ in $\R^4=\C\times\C$ given by
$$
\g(t)= (re^{ikt}, R e^{imt}), t\in[0;2\pi] ,\ k,m\in\mathbb N.
$$
This is a toric knot in a 3-sphere.
\end{example}

We have an immediate
\begin{proposition}
1)	The parameter $t$ on $\g$ is proportional to the arc-length.\\
2) The function $L(x,y)$ depends only on the difference $y-x$.
\end{proposition}

\begin{corollary} \label{cor:int}
1)	For any chord $[\g(x),\g(y)]$, the angles $\alpha, \beta$ (see formulas (\ref{first})) formed by the chord with $\g$ are equal.\\
2) The wire billiard for $\g$ is totally integrable at least near the boundary. The angle $\alpha$ is a conserved quantity.
\end{corollary}

Let us note that the curve $\g$ does not necessarily meets the niceness condition 3 of Definition \ref{nice} and hence not always determines a twist map of the cylinder.
However, in many cases it does.  We demonstrate this with the following important example which can be easily generalized.

\begin{theorem} \label{thm:niceint}
Consider the curve $\gamma$ in $\R^4=\C\times\C$ given by
\begin{equation}\label{eq:g}
\g(t)= (e^{it}, \varepsilon e^{imt}), t\in[0;2\pi] ,\ 1<m\in\Z
\end{equation}
Then, for sufficiently small $\varepsilon$, the curve $\g$ is nice, and the wire billiard is totally integrable: the angle $\alpha$ is conserved by the billiard map.
\end{theorem}

\proof
This follows immediately from Theorem \ref{thm:open} and  Corollary \ref{cor:int}.
\proofend

Consider the curve $\g$ given by (\ref{eq:g}) for $m\ge 3$. We fix a number $d>0$ and consider the ruled surface $S_{d}$ determined by the  invariant curve corresponding to $d$:

$$
S_d(t,s)=\g(t)+s\cdot r(t),
$$
where $r(t)$ is a unit vector given by
\begin{equation}\label{eq:L}
r(t)=\frac{R(t)}{L}, R(t)={\g(t+d)-\gamma(t)} , \ L={|\g(t+d)-\gamma(t)|}.
\end{equation}
Note that, by the definition of $\gamma$, the length $L$ does not depend on $t$.
Let $\delta_d$ be the striction curve of the ruled surface $S_d$ (see Theorem \ref{thm:striction}).

\begin{theorem}\label{thm:caustic}
The striction curve $\delta_d$ is a genuine caustic (the envelope of the family of chords) if and only if the relation
\begin{equation} \label{eq:Gutkin}
\tan\left(m\frac{d}{2}\right)=m\tan \left(\frac{d}{2}\right)
\end{equation}
holds.
\end{theorem}

\proof
1. Let us check first that the striction point of every chord is its midpoint.
We use the formula for the striction point
$$
\delta_d(t)=\g(t)+s^*\cdot r(t), \ s^*=-\frac{\dot\g\cdot\dot r}{\dot r^2},
$$
which holds true since $r$ is  unit  and $t$ is proportional to the arc-length parameter.

Since $L$ does not depend on $t$, we pass to $R$ using (\ref{eq:L}):
\begin{equation}\label{eq:sstar}
s^*=-\frac{\dot\g\cdot\dot R}{\dot R^2}L.
\end{equation}
We have:
\begin{equation*}
\begin{split}
R(t)=(e^{i(t+d)}, \varepsilon e^{im(t+d)})-(e^{it}, \varepsilon e^{imt})=&(e^{it}(e^{id}-1), \varepsilon e^{imt}(e^{imd}-1));\\
\dot R=(ie^{it}(e^{id}-1), im\varepsilon e^{imt}(e^{imd}-1))&;\
\dot\g=(ie^{it}, im\varepsilon e^{imt}).
\end{split}
\end{equation*}
Note that it is enough to compute $s^*$ for $t=0$ since it does not depend on $t$:
$$
\dot\g(0)=(i, i\varepsilon m),\quad \dot R(0)=(i(e^{id}-1)), i\varepsilon m (e^{imd}-1)).
$$
Also
$$
L=2\sqrt{\sin^2\frac{d}{2}+\varepsilon^2\sin^2\frac{md}{2}}.
$$
Passing from $\C\times\C$ to  $\R^4$ yields
$$
\dot\g(0)=(0,1,0,\varepsilon m),
$$
$$
\dot R(0)=(-\sin d,\cos d-1, -\varepsilon m\sin md, \varepsilon m(\cos md -1)).
$$
From these formulas we get:
$$
\dot\g\cdot \dot R=\cos d-1 +\varepsilon^2m^2(\cos md-1),\ \dot R^2=2((1-\cos d)+\varepsilon^2m^2(1-\cos md)).
$$
Therefore, by (\ref{eq:sstar}), we have
$$
s^*=L/2.
$$

2. Using the first step,  we write the striction curve in the form:
$$
\delta_d(t)=\frac{1}{2}\left((e^{it}+e^{i(t+d)}),\varepsilon(e^{imt}+e^{im(t+d)})\right).
$$
Thus
$$
\dot\delta_d(t)=\frac{1}{2}\left((ie^{it}(e^{id}+1), im\varepsilon e^{imt}(e^{imd}+1))\right).
$$
Also we have already computed
$$
R=(e^{it}(e^{id}-1), \varepsilon e^{imt}(e^{imd}-1)).
$$
Comparing last two formulas we get
$$
R\parallel \dot \delta_d\quad \Leftrightarrow\quad \frac{1-e^{id}}{1+e^{id}}=\frac{1}{m}\cdot\frac{1-e^{imd}}{1+e^{imd}}\quad \Leftrightarrow\quad \tan\left(m\frac{d}{2}\right)=m\tan \left(\frac{d}{2}\right).
$$
\proofend

\begin{remark} \label{rm:Gut}
{\rm
Equation (\ref{eq:Gutkin}) appeared earlier in the study of billiards and of flotation problems \cite{Gu1,Gu2,Gu3}. In particular, consider a convex planar billiard that possesses an invariant curve that consists of chord that make a constant angle $0<\alpha<\pi/2$ with the boundary of the billiard table at both endpoints. If the table is not circular, then the angle $\alpha$  satisfies the equation
$$
\tan (m\alpha) = m \tan\alpha
$$
for some integer $m \ge 2$.
}
\end{remark}

\section{Ellipsoids} \label{sect:ell}

Let $M$ be the ellipsoid in $\R^n$ given by $\{Ax \cdot x =1\}$, where $A$ is a positive self-adjoint operator. Let $\g(t)$ be an arc length parameterized geodesic on $M$, and let $\g(0)=x, \g'(0)=v$.

\begin{lemma} \label{curv}
The curvature $k$ of the curve $\g(t)$ at point $\g(0)$ is given by the formula
$$
k=\frac{Av\cdot v}{|Ax|}.
$$
\end{lemma}

\proof
The vector $Ax$ is a normal to $M$ at point $x$, and $N=Ax/|Ax|$ is the outward unit normal. One has $A\g\cdot \g'=0$ and, by differentiating, $A\g'\cdot \g' + A\g\cdot \g'' =0$. But $\g''=-k N$, and evaluating at $t=0$, one gets $k |Ax| =  Av\cdot v$, as needed.
\proofend

Let $M_\lambda$ be a family of confocal ellipsoids
$$
\sum_{i=1}^n \frac{x_i^2}{a_i^2+\lambda} =1
$$
with distinct semi-axes $a_i$. Then each ellipsoid in this family serves as a caustic for the billiard inside any other ellipsoid from this family that contains it.

 \begin{figure}[hbtp]
\centering
\includegraphics[width=2.5in]{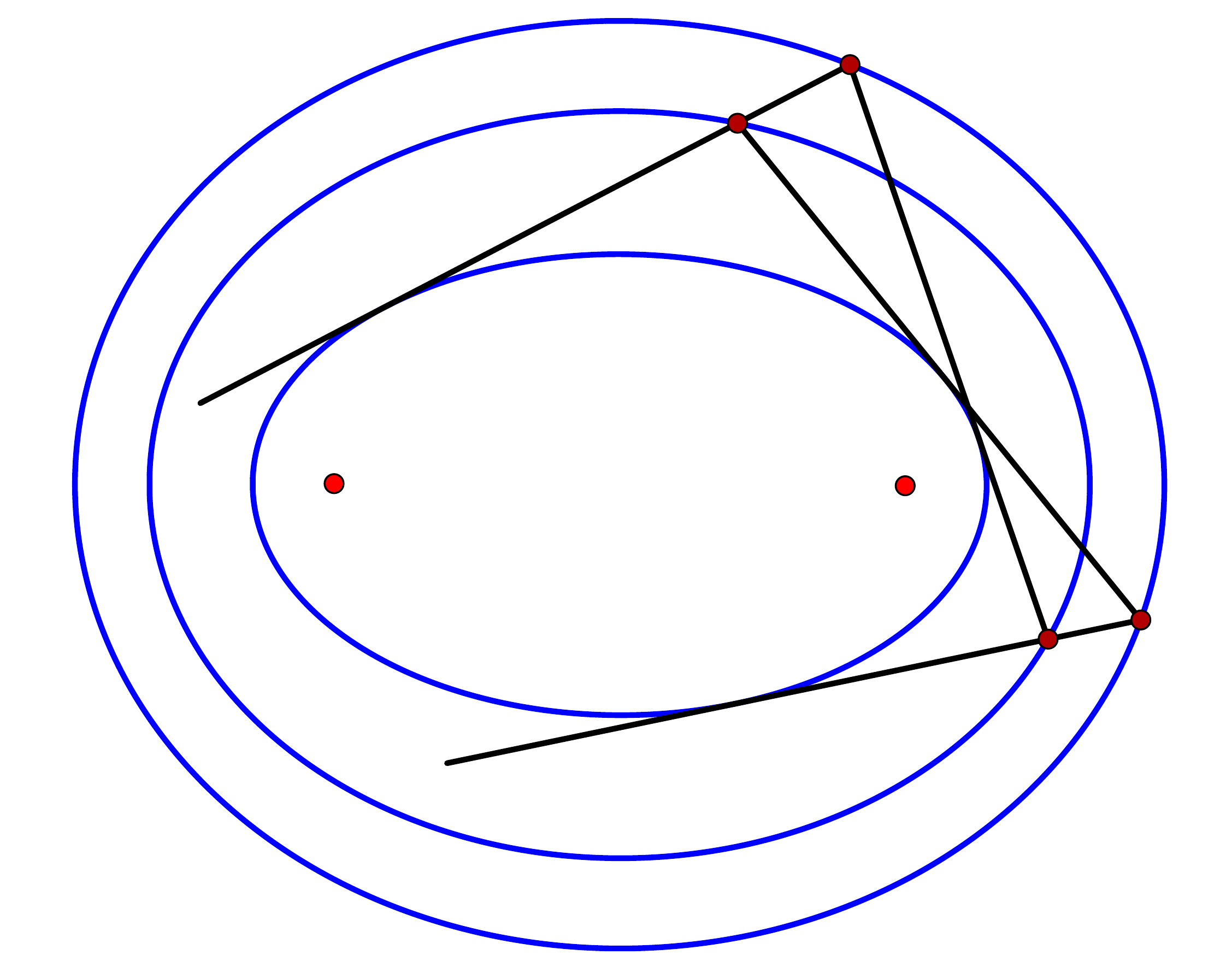}
\caption{Commuting billiard transformations on confocal ellipses.}
\label{commute}
\end{figure}

The following properties hold (see, e.g., \cite{Tab95}):
\begin{enumerate}
\item The billiard reflection in an arbitrary hypersurface, considered as a transformation of the space of oriented lines, is symplectic.
\item The billiard reflections in confocal ellipsoids commute, see Fig. \ref{commute}. This is because these transformations share integrals that Poisson commute with respect to the same symplectic structure $\Omega$.
\item The set of oriented lines, tangent to a hypersurface $M$, is a hypersurface in the symplectic space of oriented lines. The characteristic curves of this hypersurface consist of the lines tangent to a fixed geodesic of $M$.
\item As a consequence, the billiard reflection in a confocal ellipsoid, say $M_\lambda, \lambda>0$, takes the geodesics of the ellipsoid $M_0$ to the geodesics.
\end{enumerate}

\begin{figure}[hbtp]
\centering
\includegraphics[width=3in]{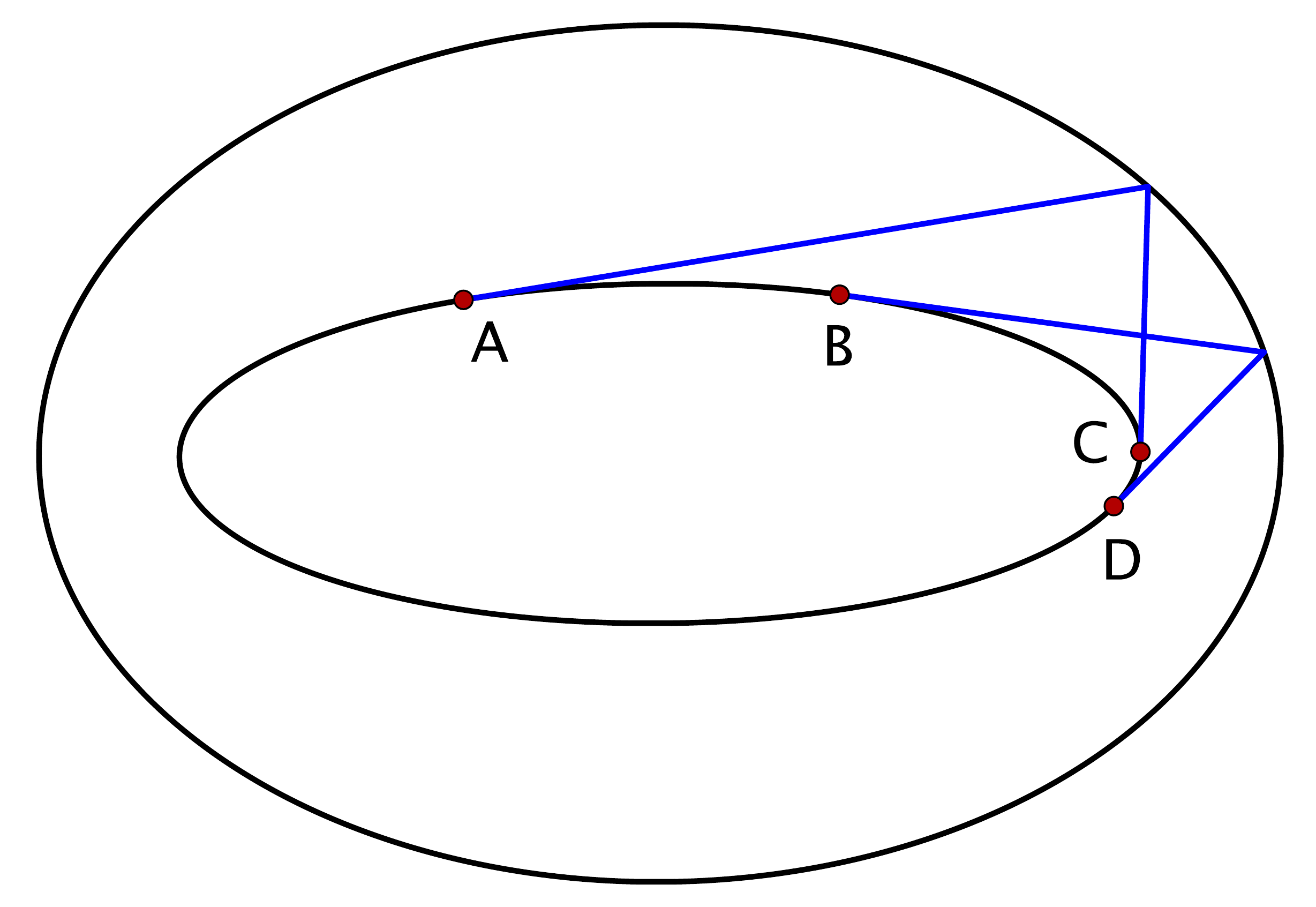}
\caption{The arc $AB$ is visibly longer than the arc $CD$. They have the same length with respect to the measure $k^{2/3} dt$.}
\label{length}
\end{figure}

However, this transformation of geodesics does not preserve the arc length parameterization, that is, the billiard reflection in a confocal ellipsoid $M_\lambda$ does not commute with the geodesic flow on $M_0$, see Fig. \ref{length}. In dimension 2, where confocal ellipses are given by the string construction, the parameter $x$, that is invariant under the billiard reflection in confocal ellipses, is given by $dx=k^{2/3} dt$, see \cite{Gl19,Por}.

Let $\xi$ be a reparameterization of the geodesic flow on the ellipsoid with the speed of the geodesic through point $x$ in the direction of unit vector $v$ is $k^{-2/3}$, where the curvature $k$ is given by Lemma \ref{curv}. We consider $\xi$ as a vector field on the set of oriented lines tangent to the ellipsoid.

\begin{theorem} \label{thm:repar}
The billiard reflections in the confocal ellipsoids $M_\lambda$ commute with the flow of the vector field $\xi$ on $M_0$.
\end{theorem}

\proof
Let $\eps$ be small. The billiard reflections in $M_\lambda$ and $M_\eps$ commute. The billiard reflection in $M_\eps$ is approximated by the reparameterized geodesic flow $\xi$, and taking the limit $\eps\to 0$ yields the result.
\proofend

\section{Open problems} \label{sec:problems}

\subsection {Sub-Riemannian approach}
Baryshnikov and Zharnitsky \cite{BZ} and, independently, Landsberg \cite{La} applied methods of sub-Riemannian geometry in the study of billiards; later these ideas were used, by a number of authors, in the study of outer billiards and of symplectic billiards. 

In a nutshell, the  idea is as follows. Let $P$ be a planar $k$-gon which is an $k$-periodic billiard orbit. Then one knows the tangent directions to the boundary of the billiard table at the vertices of $P$: they are orthogonal to the bisectors of the respective angles of the polygon. This defines an $k$-dimensional distribution on the $2k$-dimensional space of $k$-gons, called the {\it Birkhoff distribution}. The Birkhoff distribution is tangent to the level hypersurfaces of the perimeter function and is completely non-integrable therein.

If a billiard has an invariant circle consisting of $k$-periodic points, then the Birkhoff distribution on $k$-gons has a closed horizontal curve. This  makes it possible to construct billiards possessing invariant circles consisting of periodic points.

One can define an analog of Birkhoff distribution on the space of $k$-gons in $\R^n$. The space of $k$-gons is $(kn)$-dimensional, and the distribution has codimension $k$: the non-holonomic constraint on each vertex of the polygon is that its velocity is orthogonal to the bisector of the polygon at this vertex. As in the plane, this distribution is tangent to the level hypersurfaces of the perimeter function.
We expect this distribution to have properties and applications similar to the planar case.

We propose to study a more restrictive non-holonomic constraint on polygons in $\R^n$ that reduce the number of degrees of freedom of each vertex to one. Let $P=(X_1,\ldots,X_k)$ be  a $k$-gon, and let $V_i$ be the velocity of vertex $X_i$. In addition to $V_i$ being orthogonal to the bisector of $P$ at $X_i$, we require the three vectors $V_i, X_{i+1}-X_i, V_{i+1}$ to be coplanar for $i=1,\ldots,k$.

If a wire billiard has a caustic (and not just a ruled surface) corresponding to an invariant circle consisting of $k$-periodic points, then one has a closed curve in the space of $k$-gons satisfying these non-holonomic constraints. It would be interesting to develop these ideas and to construct such wire billiards.

\subsection{Ivrii conjecture}
A version of the Ivrii conjecture asserts that, for a billiard with a smooth boundary in the Euclidean plane, the set of periodic  trajectories has zero measure. This is obvious for period two, and for period three there are several proofs \cite{BZ,La,Ry,St,Vo,Wo}. Currently, the strongest result concerns period four \cite{GlK}.

We ask whether the Ivrii conjecture holds for wire billiards. If the map is multi-valued, one may have counterexamples, see Example \ref{ex:perp}. If the wire is nice, the case of period two still holds trivially, and we wonder whether the next case, period three, holds as well.

\subsection{Integrable cases. Hopf rigidity}

In the planar case, two examples of integrable billiards are known: a circle and an ellipse. The billiard in a circle is integrable in a stronger sense: the phase cylinder is foliated by invariant circles. The converse statement also holds, see \cite{Bia1}.

The billiard inside an ellipse is integrable with a more complicated phase portrait: only a part of the phase cylinder, containing both boundaries, is foliated by invariant circles; the respective caustics are confocal ellipses. The  conjecture, attributed to Birkhoff, asserts that if a neighborhood of the boundary of a convex planar billiard is foliated by caustics, then the billiard is elliptic. In spite of the efforts of many a mathematician, this conjecture remains open.  For a survey of the recent progress in this conjecture see \cite{KS}.

We conjecture that an analog of Bialy's theorem \cite{Bia1} holds for wire billiards: if the phase cylinder is foliated by invariant circles, then the wire is an orbit of a 1-parameter subgroup of the orthogonal group, cf. Theorem \ref{thm:niceint}. Notice that this example exists only for dimensions higher than 3, so no integrable examples are known in $\R^3$.
 
 One also expects to have a version of another theorem by Bialy \cite{Bia2} that estimates from above the measure of the part of the phase cylinder consisting of the minimal billiard orbits; the estimate is in terms of the length and curvature of the billiard curve and the area bounded by it.

We also wonder whether there exist non-planar analog of ellipses or other integrable billiards with more complicated phase portraits.

It would be interesting to develop an algebraic approach to the search of integrable wire billiards. In this approach one restricts to the case of algebraic curve $\g$ and looks for an additional  (rational?) function defined on the variety of all chords of $\g$,  which is invariant under wire billiard law.
This would give a further extension of the algebraic approach to billiards, see \cite{BM} for a recent survey. 

\subsection{Periodic orbits}
If the wire billiard correspondence is multi-valued, it still makes sense to estimate from below the number of its periodic orbits. Namely, given a closed curve $\g \subset \R^n$, one wants to estimate from below the number of inscribed $k$-gons of extremal length. A more general problem, when the curve $\gamma$ is replaced by an immersed submanifold, was studied in \cite{Du}; see also \cite{Pu} for period two, that is, binormal chords of a submanifold.

It is interesting to study how complexity of the wire affects the number of periodic orbits. For example, in dimension 3, one expects a complex knot to have many such inscribed polygons. Let us mention that the binormal chords of knots appear, as generators, in the recently developed knot contact homology theory.

For polygons inscribed into a closed curve, one has the notion of the rotation number, defined in the same way as in the planar case.  We conjecture that if the curvature of the curve $\gamma$ does not vanish then, for every period $q\ge 2$ and every rotation number $p \le q/2$, there exist at least two $q$-periodic wire orbits with rotation number $p$. For smooth strictly convex planar billiards, this is the celebrated theorem of Birkhoff. We believe that a modification of the argument in \cite{FaTa} should apply in this situation. 

\subsection{Other dimensions and other geometries}
Finally, it is interesting to study a version of wire billiards with a wire replaced by a submanifold of higher dimension. In particular, what are analogs of the niceness conditions, sufficient for the correspondence to be single-valued?

Another natural problem is to extend wire billiards to non-Euclidean set-up, for example, to the spherical and hyperbolic geometries.

\end{document}